\documentclass[12pt]{amsart}
\usepackage{amsmath,amssymb, graphics, amscd,latexsym}
\usepackage{epsfig}

\makeatletter

\newtheorem{Theorem}{Theorem}[section]
\newtheorem{Lemma}[Theorem]{Lemma}

\newtheorem{Proposition}[Theorem]{Proposition}
\newtheorem{Definition}[Theorem]{Definition}

\theoremstyle{remark}

\newtheorem{Example}[Theorem]{Example}

\usepackage{pb-diagram}
\begin{document}

\title{The Geometry of the Hilton Splitting}

\author{Jianhua Wang}

\begin{abstract}  One of the important theorems in homotopy theory
    is the Hilton Splitting Theorem which states: there is an
    isomorphism $H=\oplus_{\gamma\in\Gamma} H_{\gamma}$ from
    the $m$-th homotopy group of the wedge of a number of spheres
    to the direct sum of the $m$-th homotopy groups of some spheres,
    see [Hi]. In this paper we will construct geometrically all
    Hilton homomorphisms $H_{\gamma}$ and prove a family of 
    sharper symmetry relations between linking coefficients which
    desuspend and generalize the relations of Kervaire [Ke], Haefliger
    and Steer [Ha,St]. 
\vskip 4pt
\leftline {{\bf Keywords:} basic Whitehead products, Hilton 
    homomorphisms, framed links, linking coefficients.}
\vskip 4pt
\leftline {{\bf AMS Classification:} primary 55Q20, 
    57Q45; secondary 55Q25, 57R15, 57R19.}
\end{abstract}

\maketitle

\section{introduction} 
\vskip 3pt
         Let $X'$ be a connected smooth manifold without boundary
      and let $X=X'\times \Bbb R$ be of dimension $m\ge 2$. We
      denote by $k_1, k_2,\cdots, k_r$ some natural numbers greater
      than $1$. A framing of a $k$-codimensional submanifold in
      $X\times \Bbb R$ is a trivilization of the normal vector bundle,
      or equivalently, an ordered set $(u_1, u_2, \cdots, u_k)$ of
      $k$ linearly independent normal vector fields. 
      A $(k_1, k_2,\cdots, k_r)$-link is a disjoint union 
      $M_1\sqcup M_2\sqcup\cdots\sqcup M_r\subset X\times\Bbb R$ of 
      closed, framed submanifolds of codimensions  
      $k_1, k_2,\cdots, k_r$. We denote the bordism group of such
      framed links by $FL_X^{k_1, k_2,\cdots, k_r}$ which is 
      isomorphic to the homotopy group 
      $[\Sigma X_C, \vee_{i=1}^r S^{k_i}]$
      via the Pontryagin-Thom construction, where $\Sigma X_C$
      is the suspension of the one point compactification of $X$.
      Denote by $\iota_i:S^{k_i}\hookrightarrow \vee_{i=1}^r S^{k_i}$
      the inclusions and by $\Gamma$ a system of basic Whitehead
      products in $\iota_1, \iota_2, \cdots, \iota_r$. So we have
      the Hilton isomorphism (generalized by Milnor, see  [Mi])
             $$ H=\oplus_{\gamma\in\Gamma} H_{\gamma}:
                FL_X^{k_1, k_2,\cdots, k_r}\longrightarrow
                \oplus_{\gamma\in\Gamma} FL_X^{q(\gamma)+1}, $$
       where $q(\gamma)$ is the height of $\gamma$, see [Hi].           
       Let $p_{\gamma}$ be the projection from the direct sum onto
       the factor $FL_X^{q(\gamma)+1}$ corresponding to $\gamma$,
       then $H_{\gamma}=p_{\gamma}\circ H$. The homomorphism $\gamma_*$,
       induced by $\gamma$, embeds this factor into 
       $FL_X^{k_1, k_2,\cdots, k_r}$. Clearly, the Hilton homomorphisms
       $H_{\gamma}$ are characterized by

            (a) $H_{\gamma}\circ\gamma_*=id$, for $\gamma\in\Gamma$;

            (b) $H_{\gamma}\circ\gamma'_*=0$, for 
                $\gamma,\gamma'\in\Gamma$ and $\gamma'\not =\gamma$.

           As the main result of this paper we will prove a family
        of symmetry relations between linking coefficients and construct 
        geometrically all Hilton homomorphisms $H_{\gamma}$.
        It is a classical subject to construct or interpret homotopical 
        invariants by means of differential topology, for example by
        using the well known Pontryagin-Thom construction and 
        transversal intersections of submanifolds.
        Hopf invariants and Hilton homomorphisms are of particular
        interest. For example, if the $H_{\gamma}$'s are already
        constructed, then for a given element in the homotopy group
        of the wedge of spheres we can compute its Hilton splitting.
        Kervaire [Ke] gave a geometrical description of the
        stable Whitehead-Hopf invariant and proved 
        a symmetry relation between linking coefficients in
        high dimensions. Haefliger 
        and Steer [Ha,St] constructed the suspension of $H_{\gamma}$
        corresponding to $\gamma=[\iota_1,\iota_2]$ and obtained a 
        further symmetry relation between linking coefficients. 
        Boardman and Steer [Bo,St] defined the Hopf ladder and
        presented a geometrical discussion. Koschorke and Sanderson
        applied immersion theory to this topic in [K,S 1] and [K,S 2].

           This work is also strongly motivated by the close connection
        to homotopy theory of link maps. For example, Koschorke [Ko 1] and
        [Ko 3] generalized the $\mu$-invariants of Milnor by using his
        geometrical interpretation of some stable Hilton homomorphisms
        which in fact makes some computations possible. Hilton splitting
        also played an important role in the study of the different 
        homotopy behaviour of link maps in $S^m$ and $\Bbb R^m$,
        see [Ka 1] and [Ka 2] of Kaiser.

           In the case $X\times \Bbb R\cong\Bbb R^3$ and 
        $k_1=k_2=\cdots=k_r=2$ it is well known that the elements in
        the factors $\pi_3(S^3)\cong\Bbb Z$ in the Hilton splitting of
        $\pi_3(\vee_{i=1}^r S_i^2)$ can be interpreted as linking
        numbers. Sanderson [Sa 2] gave a geometrical isomorphism
        from $\pi_4(S^2\vee S^2)$ to $\Bbb Z_2^3\oplus\Bbb Z^2$
        by using intersections with Seifert surfaces. This isomorphism
        takes the form of the Hilton splitting but is different from it,
        see the author's dissertation [Wa]. In general cases the 
        geometry of the Hilton splitting is unknown up to date, because 
        of its complicated algebraic topological nature.

          This paper is organized as follows. We introduce in \S 2
        a new construction, call it the $\tau$-construction, and establish
        its basic properties. Our $\tau$-construction desuspends the one
        of Haefliger and Steer in [Ha,St]. As an application of our
        construction we prove a family of sharper symmetry relations
        between linking coefficients in \S 3. The $\tau$-construction 
        leads to the
        definition of the $\tau$-reduction in \S 4, all Hilton
        homomorphisms are geometrically constructed there by means of
        the $\tau$-reductions. We work in the category of smooth manifolds.

           We have extracted the materials in \S 2 and \S 3 from the
        author's dissertation [Wa]. So it is a great pleasure to express
        my gratitudes to my supervisor Professor U.\ Koschorke as well as
        Professor U.\ Kaiser for many helpful discussions and encouragements.
        I am grateful to Professor M.\ Heusener for nice talks. Thanks
        also to Dr. Pho who helped me use {\it xfig}.

\vskip 1cm
\section{the $\tau$-construction}
\vskip 3pt

           Let $M_1\sqcup M_2\sqcup\cdots\sqcup M_r\subset X\times\Bbb R$
        be a $(k_1, k_2,\cdots, k_r)$-link and $1\le i\not= j\le r$. 
        We construct now a framed submanifold 
        $Z=\tau (M_j,M_i)\subset X\times\Bbb R$ as follows.
        Let $W_j=M_j\times [0,1]$, and let 
        $W_i\subset X\times\Bbb R\times [0,1]$  
        be a framed bordism of $M_i$ such that $M_j\times\{1\}$ and
        $W_i\pitchfork X\times\Bbb R\times\{1\}$ are separated
        by some $X_t=X\times\{t\}\times\{1\}$, see Fig.1.

          Denote the naturally framed intersection of $W_i$ and $W_j$
        by $\bar Z$ and let $u_{k_j}$ be the last vector field in the
        framing of $M_j$. For $\varepsilon >0$ small enough we deform
        first $W_j$ to $M_j\times [0,\varepsilon]$ and then rotate at every
        point $x\in M_j$ the interval $[0,\varepsilon]$ to $u_{k_j}$ through
        the angle $\pi/2$. By doing this we have isotoped $\bar Z$ to
        a submanifold $Z=\tau (M_j,M_i)\subset X\times\Bbb R\times\{0\}$
        which is naturally framed, because the isotopy induces a homotopy 
        of the normal vector bundles and during the isotopy $u_{k_j}$
        is deformed to the negative direction of the interval $[0,1]$. 
        See Fig.1 again. We call this construction of
        $Z=\tau (M_j,M_i)\subset X\times\Bbb R\times\{0\}$ the 
        $\tau$-construction, which desuspends the one of
        Haefliger and Steer [Ha,St].

           If one changes the roles of $M_i$ and $M_j$, namely 
        takes $W_i$ to be the cylinder $M_i\times [0,1]$ and
        takes $W_j$ to be a framed bordism, then one will get another
        framed submanifold 
        $\tau (M_i,M_j)\subset X\times\Bbb R\times \{0\}$. Denote by
        $\tau [*,*]$ the framed bordism class. According to [Ha,St] it
        holds $E\tau [M_j,M_i]=E\tau [M_i,M_j]$ up to sign, where $E$
        denotes the suspension homomorphism. But, as we will see later
        in this section, $\tau [M_j,M_i]\not=\tau [M_i,M_j]$, even if up
        to involution (namely an automorphism $u$ of the
        target group with the property $u\circ u=id$).

\begin{figure}[htb]
\setlength{\unitlength}{1bp}
\begin{picture}(255,255)(-40,0)
\epsfxsize=12cm
\put(-80,120){\epsfbox{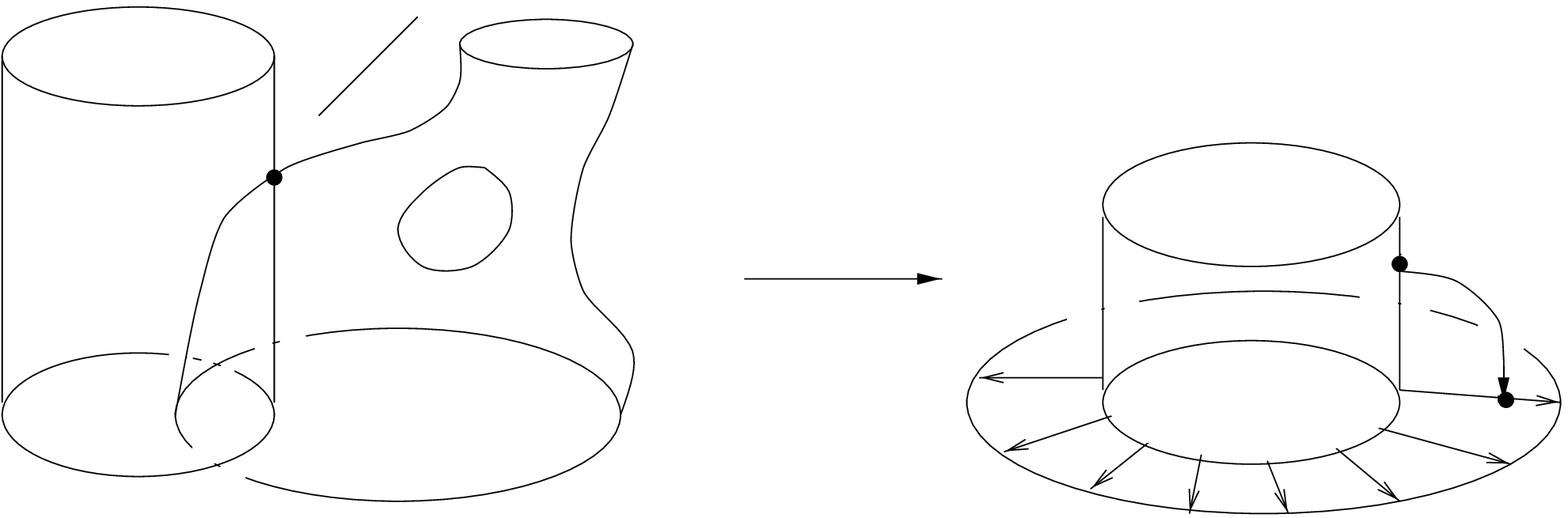}}
\put(-10,230){$X_t$}
\put(-60,180){$W_j$}
\put(-10,170){$W_i$}
\put(-20,180){$\bar Z$}
\put(260,135){$u_{k_j}$}
\put(230,150){$ Z$}
\end{picture}
\vspace{-4cm}
\caption{}
%\label{}
\end{figure}

           Because $\tau (M_j,M_i)$ (or $\tau (M_i,M_j)$) lies in a
        small neighbourhood of $M_j$ (or $M_i$), it is disjoint from
        all components of the original link. This interesting fact 
        makes it possible to get a new link 
             $$M_1\sqcup M_2\sqcup\cdots\sqcup M_r\sqcup\tau (M_j,M_i)
                \subset X\times\Bbb R$$
        from the old one. Denote by $[ \ ]$ the bordism class of a
        framed submanifold or link.
\begin{Theorem}
The assignment
        $$ M_1\sqcup M_2\sqcup\cdots\sqcup M_r\to
        M_1\sqcup M_2\sqcup\cdots\sqcup M_r\sqcup\tau (M_j,M_i) $$
        gives a well defined injective homomorphism $\tau_{ji}$
\[
\begin{diagram}
         \node{FL_X^{k_1,k_2,\cdots,k_r}}
            \arrow{e,t}{\tau_{ji}}
            \arrow{se,b}{\tau_{ji}^p}
         \node{FL_X^{k_1,k_2,\cdots,k_r,k_{r+1}}}
              \arrow{s,r}{proj}                       
\\
         \node[2]{FL_X^{k_{r+1}}}
\end{diagram}
\]  
        where $k_{r+1}=k_i+k_j-1$. In particular, 
        $\tau_{ji}^p=proj\circ\tau_{ji}$ is a well defined invariant.
\end{Theorem}
\begin{proof} \ Let $I_0=[0,1]$, $I_1=[1,2]$, $I_2=[-1,0]$ and let 
        $N_1\sqcup N_2\sqcup\cdots\sqcup N_r\subset X\times\Bbb R\times I_0$
        be a framed bordism from $M_1'\sqcup M_2'\sqcup\cdots\sqcup M_r'$
        to $M_1\sqcup M_2\sqcup\cdots\sqcup M_r$. We perform
        $\tau (M_j,M_i)$ in $X\times\Bbb R\times I_1$ using the cylinder
        $W_j$ and similarly we perform $\tau (M'_j,M'_i)$ in 
        $X\times\Bbb R\times I_2$ using the cylinder $W'_j$ (the negative
        orientation of $I_2=[-1,0]$ is used), see Fig.2.
\begin{figure}[htb]
\setlength{\unitlength}{1bp}
\begin{picture}(255,255)(-40,0)
\epsfxsize=15cm
\put(-5,-75){\epsfbox{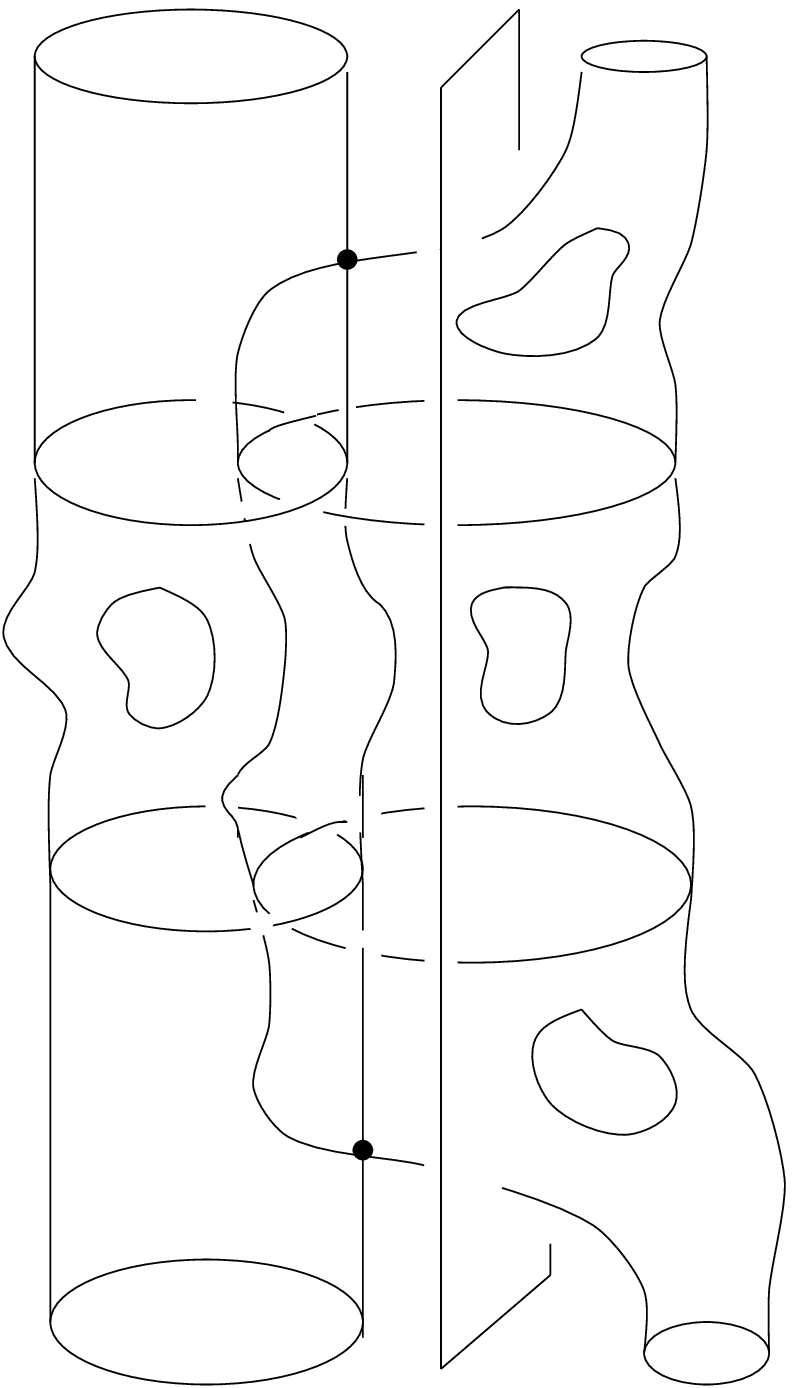}}
\put(75,220){$X'_t$}
\put(20,160){$W_j$}
\put(120,185){$W_i$}
\put(55,170){$\bar Z$}
\put(20,55){$N_j$}
\put(100,55){$N_i$}
\put(20,-20){$W'_j$}
\put(120,-30){$W'_i$}
\put(58,-15){$\bar Z'$}
\put(-25,100){$V_j$}
\put(150,100){$V_i$}
\put(165,-70){$A$}
\put(150,205){$B$}
\end{picture}
\vspace{3.cm}
\caption{}
%\label{}
\end{figure}

           Consider now the framed submanifolds $V_j=W'_j\cup N_j\cup W_j$
        and  $V_i=W'_i\cup N_i\cup W_i$. For $t\in\Bbb R$ let 
        $\Bbb R_t^{-}=\{x\in\Bbb R | x<t\}$ and similarly $\Bbb R_t^{+}$.
        We may assume that there is some $t\in\Bbb R$ such that
\begin{align*}
              V_j & \subset   X\times\Bbb R_t^-\times [-1,2] \\
              A   & =         W'_i\pitchfork X\times\Bbb R\times\{-1\}
                              \subset 
                               X\times\Bbb R_t^+\times\{-1\}, \\
              B   & =          
                               W_i\pitchfork X\times\Bbb R\times\{2\}
                               \subset
                               X\times\Bbb R_t^+\times\{2\}.   
\end{align*}
        This means that $V_j$ and the boundary of $V_i$, namely the
        union of $A$ and $B$, are separated by
        $X'_t=X\times\{t\}\times [-1,2]$, see Fig.2 again. 
        Because we are working in $X\times\Bbb R=X'\times\Bbb R^2$,
        this can always be satisfied by isotoping
        $V_i$ without changing the framed intersection 
            $$ V_j\pitchfork V_i=\bar Z\sqcup \bar Z'
               =(W_j\pitchfork W_i)\sqcup (W'_j\pitchfork W'_i). $$

\begin{figure}[htb]
\setlength{\unitlength}{1bp}
\begin{picture}(255,255)(-40,0)
\epsfxsize=10cm
\put(-45,50){\epsfbox{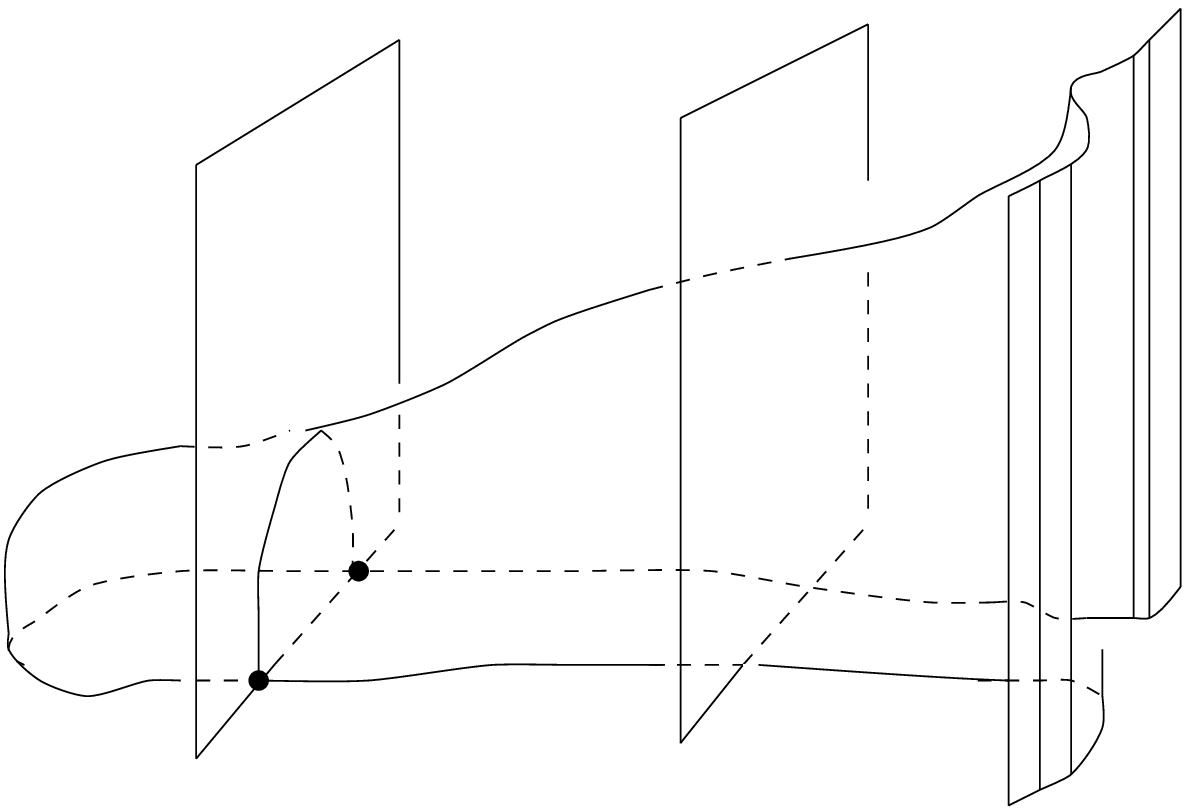}}
\put(115,200){$X''_t$}
\put(20,180){$Q_j$}
\put(10,125){$\bar Q$}
\put(3,90){$\bar Z'$}
\put(42,93){$\bar Z$}
\put(90,150){$Q_i$}
\put(165,150){$A\times I$}
\put(240,160){$B\times I$}
\put(170,70){$V_i$}
\put(190,215){$V'_i$}
\end{picture}
\vspace{-2.cm}
\caption{}
%\label{}
\end{figure}
                            
         Let $I=[0,1], Q_j=V_j\times I$ and $Q_i\subset X\times\Bbb R
      \times [-1,2]\times I$ be a framed bordism of $V_i$ such that
            $$ \partial Q_i=V_i\cup A\times I\cup B\times I
                            \cup V'_i  $$
      and such that $Q_j, V'_i$  are separated by 
      $X_t''=X\times \{t\}\times [-1,2]\times I$, where $V'_i$ is the
      boundary part of $Q_i$ lying in 
      $X\times \Bbb R\times [-1,2]\times\{1\}$, see Fig.3. 
      Such a manifold $Q_i$ always exists.
        
         Let $\bar Q=Q_j\pitchfork Q_i$ be the naturally framed 
      intersection, its boundary is $\partial\bar Q=\bar Z\sqcup\bar Z'$.
      Just because of 
         $$ Q_j=(M'_j\times [-1,0]\times I)\cup (N_j\times I)
                 \cup (M_j\times [1,2]\times I)   $$
      there is an isotopy of $Q_j$ which deforms $Q_j$ to
      $N_j\times I$ and is smooth at least in a small neighbourhood of 
      $\bar Q\subset Q_j$. For example, for any $x\in M_j$ one can 
      easily isotope $\{x\}\times [1-\varepsilon,2]\times [0,1]$ to            
      $\{x\}\times [1-\varepsilon,1]\times [0,1]$ by using the trick
      in Fig.4. Note that collars of $M_j, M'_j\subset N_j$ are used to 
      perform this isotopy. 
\begin{figure}[htb]
\setlength{\unitlength}{1bp}
\begin{picture}(255,255)(-40,0)
\epsfxsize=14cm
\put(-82,150){\epsfbox{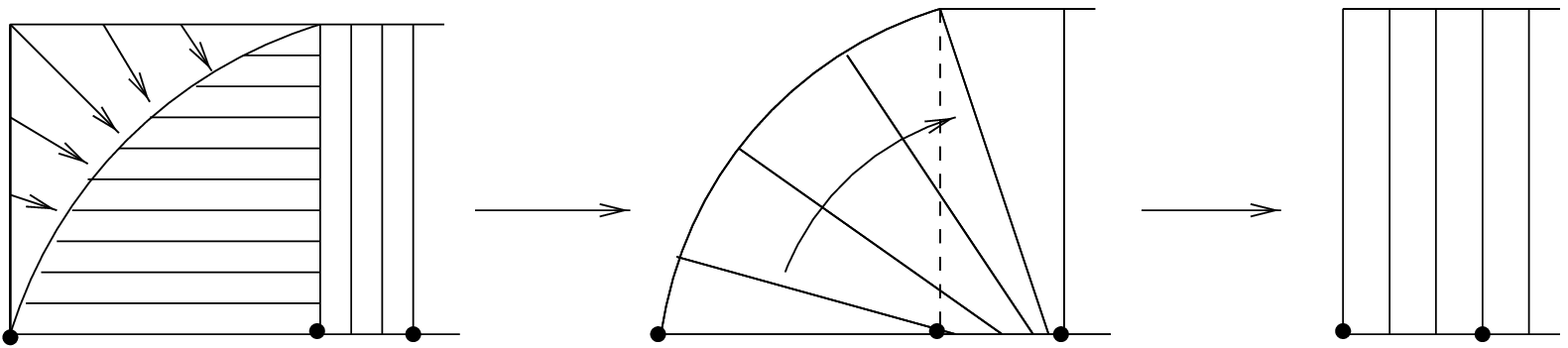}}
\put(-92,220){$1$}
\put(-92,150){$0$}
\put(-85,140){$2$}
\put(-9,140){$1$}
\put(3,140){$1-\varepsilon$}
\put(76,140){$2$}
\put(143,140){$1$}
\put(160,140){$1-\varepsilon$}
\put(241,140){$1$}
\put(265,140){$1-\varepsilon$}
\end{picture}
\vspace{-5cm}
\caption{}
%\label{}
\end{figure}

        So we can isotope $\bar Q$ smoothly to a framed submanifold
      $\bar Q'\subset N_j\times I$. Let $\tilde u_{k_j}$ be the 
      last normal vector field in the framing of $N_j$. Just like
      in the $\tau$-construction we deform $N_j\times I$ to
      $N_j\times [0,\varepsilon]$ and then rotate the positive 
      $I$-direction to $\tilde u_{k_j}$. By doing this we have
      isotoped $\bar Q'$ to a framed submanifold 
      $Q\subset X\times\Bbb R\times [0,1]\times\{0\}$. Now it is
      easy to see that $Q$ is a framed bordism between 
      $\tau (M_j, M_i)$ and $\tau (M'_j, M'_i)$, and $Q$ is also
      disjoint from $N_1\sqcup N_2\sqcup\cdots\sqcup N_r$, for it 
      lies in a small neighbourhood of $N_j$. The desired framed
      bordism is given by 
      $N_1\sqcup N_2\sqcup\cdots\sqcup N_r\sqcup Q$.

         So $\tau_{ji}$ is a well defined map. The assumption
      $X=X'\times \Bbb R$ guarantees it is also a homomorphism.
      Other assertions follow easily.
\end{proof}

       Note that the homomorphism $\tau_{ji}$ in the theorem
     is independent of the choice
     of the vector field used to rotate $M_j\times [0,\varepsilon]$,
     because one can always rotate one vector field to another.
     Theorem 2.1 implies that we can perform the $\tau$-construction
     successively to get further well-defined invariants. For example, 
     for any $1\le k\le r$ we have $\tau(M_k,\tau(M_j, M_i))$. 
     This is a very important property
     of our $\tau$-construction. To understand this, note that
     in the Haefliger-Steer construction we see only one geometrical
     obstruction (the framed intersection $\bar Z$) of a framed link of 
     two components from being the trivial link (in the sense the two
     components are not linked), in contrast the Hilton
     splitting says there are many other obstructions; our iterated 
     $\tau$-invariants are surely related to such further obstructions.  

        Let $M\subset X\times\Bbb R$ be a closed, framed submanifold.
     A suitably framed Seifert surface of $M$ is a compact, framed
     submanifold $F\subset X\times\Bbb R$ with boundary $M$ such 
     that the framing of $M$ as boundary is homotopic to the original
     framing. In this case we say $M$ is $S$-framed. Note that two
     $S$-framings of $M$ are not necessarily homotopic. 
     The first part of the
     following lemma is directly to see and the second part follows
     by a simple discussion of fibre-wise embeddings, so we omit
     the proof.

\begin{Lemma} \ (i) If $M_i$ has a suitably
    framed Seifert surface $F_i$, then the following two framed links 
      $$M_1\sqcup M_2\sqcup\cdots\sqcup M_r\sqcup\tau(M_j,M_i),\hskip 1cm
        M_1\sqcup M_2\sqcup\cdots\sqcup M_r\sqcup (M_j\pitchfork F_i)^{sh} $$
    are framed bordant at least up to involution of the framing of 
    the last component,
    where $(M_j\pitchfork F_i)^{sh}$ is a small shift of 
    $M_j\pitchfork F_i$ along the framing of $M_j$.

       (ii) \ Let $L_1\subset \Bbb R^n$ and $L_2\subset \Bbb R^{n'}$
    be $(k_1,k_2,\cdots,k_r)$-links with components $M_i$ and $M'_i$
    respectively, $Z, Z'\subset X\times \Bbb R$
    be disjoint, closed and framed submanifolds of codimensions
    $n$ or $n'$. By means of fibre-wise embeddings we 
    get a new $(k_1,k_2,\cdots,k_r)$-link $L\subset X\times\Bbb R$
    with components $\bar M_i=Z\times M_i\sqcup Z'\times M'_i$.
    If $M_i$ and $M'_i$ are framed zerobordant for some $i$, 
    then we can perform
    the $\tau$-construction so that the following holds
        $$ \tau_{ji}(L)=Z\times\tau_{ji}(L_1)
                        \sqcup Z'\times\tau_{ji}(L_2).$$
    The special case $Z'=\phi$ is also useful.
\end{Lemma} 
        
       The inclusions $\iota_i:S^{k_i}\hookrightarrow 
    \vee_{i=1}^r S^{k_i}$ and the Whitehead products 
    $[\iota_i,\iota_j]$ can be geometrically interpreted as 
    framed points or as $S$-framed Hopf links, via Pontryagin-Thom
    construction. Iteratedly we can represent every Whitehead product
    $\gamma$ in $\iota_1,\cdots,\iota_r$ by a 
    $(k_1,k_2,\cdots,k_r)$-link in $\Bbb R^{q(\gamma)+1}$,
    where $q(\gamma)$ is the height of $\gamma$. We are now ready
    to identify the homomorphism $\tau^p_{ji}$ in Theorem 2.1 with 
    the Hilton homomorphism 
    corresponding to $\gamma=[\iota_i,\iota_j]$.

\begin{Theorem} \ Let $\gamma=[\iota_i,\iota_j]\in \Gamma$ be
    a basic Whitehead product. It holds $\tau^p_{ji}=H_{\gamma}$
    up to involution.
\end{Theorem}
\begin{proof} \ We show that up to involution $\tau^p_{ji}$ 
    satisfies the properties (a) and (b) in \S 1 
    which characterize the Hilton homomorphisms. 
    For $\gamma'=\iota_k\in \Gamma$ property (b) is trivial. Assume that
    $\gamma'\in \Gamma$ is of weight $\ge 2$ and let
    $L=M_1\sqcup M_2\sqcup\cdots\sqcup M_r\subset\Bbb R^{q(\gamma')+1}$
    be the framed link representing $\gamma'$. Each component $M_i$
    is clearly framed zerobordant, so Lemma 2.2 reduces (a) and (b)
    to the following

       ($a'$) \ $\tau[M_j,M_i]=\pm 1$, 
              if $\gamma=\gamma'=[\iota_i,\iota_j],$

       ($b'$) \ $\tau[M_j,M_i]=0$, 
              if $\gamma'\not=\gamma=[\iota_i,\iota_j].$

       Let $\gamma'=[\iota_{i'},\iota_{j'}]$ be of weight $2$. If
    the pair $(i',j')\not=(i,j)$, then $M_i$ or $M_j$ is the empty,
    ($b'$) follows easily; if $(i',j')=(i,j)$, then $M_i\sqcup M_j$
    is an $S$-framed Hopf link and all other components are the empty,
    ($a'$) follows by Lemma 2.2.
 
        Let $\gamma'=[\alpha,\beta]$ be of weight $\ge 3$. We have
     the following formula
\begin{eqnarray}
            M_k=S^{l_1-1}\times M_k(\beta)\sqcup        
                S^{l_2-1}\times M_k(\alpha),
\end{eqnarray}
     where $S^{l_1-1}\sqcup S^{l_2-1}$ is an $S$-framed Hopf link 
     and $M_k(\alpha)$, $M_k(\beta)$ are the components of the
     links representing $\alpha$ and $\beta$, $1\le k\le r$.
     The weight of $\beta$ is at least $2$, so $M_i(\beta)$ is 
     framed zerobordant. If $\alpha\not=\iota_i$, then $M_i(\alpha)$
     is also framed zerobordant. We use Lemma 2.2 again and obtain         
\begin{eqnarray}
        \tau(M_j,M_i)=S^{l_1-1}\times\tau(M_j(\beta),M_i(\beta))\sqcup
                      S^{l_2-1}\times\tau(M_j(\alpha),M_i(\alpha)).
\end{eqnarray}
     By inductive assumption for $\alpha$ and $\beta$ we
     see that at least one of $\tau(M_j(\alpha),M_i(\alpha))$ and    
     $\tau(M_j(\beta),M_i(\beta))$ is framed zerobordant, say the first.
     By means of the fibre-wise embedding of the framed zerobordism
     one can easily prove that $\tau(M_j,M_i)$ is framed bordant to
     $S^{l_1-1}\times\tau(M_j(\beta),M_i(\beta))$, and which is clearly
     framed zerobordant, for $S^{l_1-1}$ is $S$-framed and therefore
     framed zerobordant. ($b'$) follows.

       Now let $\alpha=\iota_i$. It holds then
     $\gamma'=[\iota_i,[\iota_{i_1},\cdots,[\iota_{i_t},\iota_{j'}]\cdots]]$
     according to the definition of basic Whitehead products,
     where $\iota_i\ge\iota_{i_1}\ge\cdots\ge\iota_{i_t}<\iota_{j'}$.
     If $\iota_{j'}\not=\iota_j$, then $\iota_j$ does not appear in
     $\gamma'$, because $\iota_i<\iota_j$. This means $M_j=\phi $ and
     ($b'$) follows. So let $\iota_{j'}=\iota_j$. We assume now           
     $\iota_i=\iota_{i_1}=\cdots=\iota_{i_t}$, otherwise the argument
     is essentially the same. Denote by $w$ the weight of $\gamma'$.
     The link representing $\gamma'$ is given by
\begin{eqnarray}
          M_i & =  & S^{k_i-1}_{w-1}\times S^{k_i-1}_{w-2}\times
                   \cdots\times S^{k_i-1}_{2}\times S^{k_j-1}\sqcup 
                       \\ \nonumber
              &    &  S^{k_i-1}_{w-1}\times S^{k_i-1}_{w-2}\times
                    \cdots\times S^{k_i+k_j-2}\sqcup
                        \\ \nonumber
              &    &  \cdots\cdots\sqcup     \\ \nonumber
              &    &  S^{k_i-1}_{w-1}\times S^{(w-3)(k_i-1)+k_j-1}\sqcup
                         \\ \nonumber
              &    &  S^{(w-2)(k_i-1)+k_j-1}   \\ \nonumber
              & =  &  N_{i,1}\sqcup\cdots\sqcup N_{i,w-1},   \\
          M_j & =  & S^{k_i-1}_{w-1}\times S^{k_i-1}_{w-2}\times
                    \cdots\times S^{k_i-1}_{2}\times S^{k_i-1}_1, 
\end{eqnarray}
    where $S^{(k-1)(k_i-1)+k_j-1}\sqcup S^{k_i-1}_{k}$ are framed
    Hopf links, $1\le k\le w-1$. The products are given by 
    fibre-wise embeddings. All other components are the empty.

       Let $e$ be the last vector in the standard base of 
    $\Bbb R^{q(\gamma')+1}$. We may assume 
    $M_j\subset \Bbb R^{q(\gamma')}\times\{0\}\subset\Bbb R^{q(\gamma')+1}$
    and that $e$ is just the last vector field in the framing of $M_j$.
    This implies the following: the small shifts $\bar Q_1$ and 
    $\bar Q_2$ of any $Q_1, Q_2\subset M_j$ along $e$ through distances
    $d_1<d_2$ are separated by $\Bbb R^{q(\gamma')}\times\{(d_1+d_2)/2\}$.
     
       Obviously, every $N_{i,k}$ bounds some suitably framed Seifert
    surface $F_{i,k}$, $1\le k\le w-1$. In
    $\Bbb R^{q(\gamma')+1}\times [0,1]$ one can push them into different
    heights $a_1<a_2<\cdots <a_{w-1}$ (with boundaries fixed) to get a 
    framed zerobordism of $M_i$ which will be used to perform
    $\tau(M_j,M_i)$. Let $\{pt\}$ be a set consisting of a single point
    and cosider
      $$ Q_k=S^{k_i-1}_{w-1}\times\cdots\times 
             S^{k_i-1}_{k+1}\times\{pt\}
             \times S^{k_i-1}_{k-1}\times\cdots\times S^{k_i-1}_1
             \subset M_j.$$
    It is not difficult to see 
    $\tau(M_j,M_i)=\sqcup_{k=1}^{w-1}\bar Q_k$, where the framed 
    submanifolds $\bar Q_k$ are small shifts of the $Q_k$'s through
    distances $d_1<d_2<\cdots <d_{w-1}$ along $e$ (or equivalently
    along the framing). Clearly, every $\bar Q_k$ bounds a suitably
    framed Seifert surface and the discussion above shows they
    are separated. ($b'$) follows.
\end{proof}
\begin{Example} \ Our invariant $\tau_{ji}^p$ is asymmetric,
    namely $\tau_{ji}^p\not=\tau_{ij}^p$, even if up to involution.
    Let $r=2$ and consider the Whitehead products  
       $$ \gamma=[\iota_2,[\iota_1,[\iota_1,\iota_2]]], \ \ 
          \gamma_1=[\iota_1,[\iota_2,[\iota_2,\iota_1]]], \ \ 
          \gamma_2=[[\iota_1,\iota_2],[\iota_1,\iota_2]]]
                  =[\iota_1,\iota_2]\circ [\iota,\iota],     $$
    where $\iota$ is the identity of $S^{k_1+k_2-1}$. Denote the framed
    links representing $\gamma, \gamma_1, \gamma_2$ by $L, L_1$ and 
    $L_2$. From the Jacobi-identity follows $L=\pm L_1\pm L_2$.

      Let $\Gamma, \Gamma'$ be two systems of basic Whitehead products 
    in $\iota_1<\iota_2$ and $\iota_2<\iota_1$
    respectively, then $\gamma \in\Gamma$ and $\gamma'\in\Gamma'$.
    From Theorem 2.3 we have $\tau_{2,1}^p[L]=0=\tau_{1,2}^p[L_1]$ and
       $$ \tau_{1,2}^p[L]=\pm \tau_{1,2}^p[L_1]\pm\tau_{1,2}^p[L_2]
                         =\pm [\iota,\iota]\not=0 $$
    if $k_1+k_2-1\not=1,3,7$, according to a well known result of 
    G.Whitehead and F.Adams, see for example [Ad]. The statement
    is proved. This example shows that our $\tau$-construction
    catches what is lost in the Haefliger-Steer construction due
    to the suspension.
\end{Example}

       By using the $\tau$-construction successively 
    we can do the following: 
\begin{quote}
        (i) for many basic Whitehead products construct 
        homomorphisms $h'_{\gamma}$ 
        with property (a) in \S 1, but we can not guarantee the property (b), 
        so these homomorphisms may be different from the corresponding
        Hilton homomorphisms;
\newline        
        (ii) construct all the Hilton homomorphisms if the basic 
        Whitehead products of weight at least $4$ are not 
        involved in the Hilton splitting (using this we can re-prove the
        Jacobi-identity); 
\end{quote}
    for details see Kapitel 3 in the author's dissertation [Wa].

\vskip 1cm
\section{symmetry relations between linking coefficients}
\vskip 3pt

         As an application of the $\tau$-construction we prove here
    a family of symmetry relations between linking coefficients. Our
    argument is based on some beautiful ideas in [Ha,St].
  
       Let 
    $M_1\sqcup M_2\sqcup\cdots \sqcup M_r\subset \Bbb R^{m}\times\{1\}$,
    $W_j=M_j\times [0,1]$, and $W_i\subset\Bbb R^{m}\times [0,1]$      
    be as in the $\tau$-construction. If one isotopes the intersection 
    $W_i\pitchfork W_j$ into $\Bbb R^{m}\times\{0\}$ (instead of into
    $\Bbb R^{m}\times\{1\}$) and then project it to $\Bbb R^m\times\{1\}$, 
    one obtains a framed submanifold $\tau'(M_j,M_i)$. 
    By rotating the negative direction of $[0,1]$ to
    $-u_{k_j}$ one gets another framed submanifold $\tau''(M_j,M_i)$. 
    We have the following fact
\newline
    {\bf Fact:} \ The framed links
\begin{align*}
       & M_1\sqcup M_2\sqcup\cdots\sqcup M_r\sqcup\tau(M_j,M_i), \\
       & M_1\sqcup M_2\sqcup\cdots\sqcup M_r\sqcup\tau'(M_j,M_i),\\
       & M_1\sqcup M_2\sqcup\cdots\sqcup M_r\sqcup\tau''(M_j,M_i)
\end{align*} 
    are framed bordant at least up to involution of the framing
    of the last component. For the first and third links rotating $u_{k_j}$ 
    to $-u_{k_j}$ in the plane spanned by $u_{k_j}$ and any other vector
    field $u$ in the framing of $M_j$; for the first and second links 
    rotate $u_{k_j}$ through the angle $\pi$ in the plane spanned by 
    $u_{k_j}$ and $u$ with the middle point of $u_{k_j}$ fixed . 
    This fact will be used later in the proof of Theorem 3.1.

       Let $\gamma=
    [\iota_{j_t},[\cdots,[\iota_{j_2},[\iota_{j_1},\iota_1]]\cdots]]$ 
    be any Whitehead product in $\iota_1,\cdots,\iota_r$, such that
    $\iota_1$ appears in $\gamma$ exactly one time in the given position.
    Define
          $$\mu_{\gamma}^{\tau}:\pi_{*}(S^{k_1}\vee\cdots\vee S^{k_r})
            \longrightarrow \pi_{*}(S^{q(\gamma)+1})  $$
    as the framed bordism class of 
         $$Z_{\gamma}=\tau(\tau(\cdots\tau(\tau(M_1,
                      M_{j_1}),M_{j_2}),\cdots ,M_{j_{t-1}}),M_{j_t}).$$
    
      Let $L=S^{p_0}\sqcup S^{p_1}\sqcup\cdots\sqcup S^{p_r}\subset
    S^{n+1}$ be a smoothly embedded spherical link with
    $p_0, p_1,\cdots, p_r\le n-2$. From the well known $(n-1)$-homotopy
    equivalence we obtain the linking coefficients
\begin{align*}
    & \lambda_0\in\pi_{p_0}(S^{k_1}\vee S^{k_2}\vee\cdots\vee S^{k_r}), \\ 
    & \lambda_1\in\pi_{p_1}(S^{k'_1}\vee S^{k'_2}\vee\cdots\vee S^{k'_r}),
\end{align*}
    where $k_1=n-p_{1}, k'_1=n-p_0$, and $k_i=k'_i=n-p_i$ for $2\le i\le r$.
\begin{Theorem} \ Let $E$ denote the suspension homomorphism. 
    It holds 
            $$ E^{n+1-p_0}\mu_{\gamma}^{\tau}(\lambda_0)=
               \pm E^{n+1-p_1}\mu_{\gamma}^{\tau}(\lambda_1).$$
\end{Theorem}
\begin{proof} \ I. \ Let $I=I_1=[0,1]$. We assume first that following
    links are zero $h$-cobordant
\begin{align*}
         & L_0 =S^{p_0}\sqcup S^{p_2}\sqcup\cdots\sqcup S^{p_r},  \\  
         & L_1 =S^{p_1}\sqcup S^{p_2}\sqcup\cdots\sqcup S^{p_r}.
\end{align*}
    The spheres in $L_0$ bound disjoint, framed Seifert surfaces
    $V_0, V_{0,2},\cdots, V_{0,r}$; by the same token the spheres
    in $L_1$ bound disjoint, framed Seifert surfaces
    $V_1, V_{1,2},\cdots, V_{1,r}$. In addition, we can suppose that
    for $2\le i\le r$ the framed submanifolds 
         $$ V_{0,i}\times\{0\}\cup S^{p_i}\times I\cup V_{1,i}\times\{1\}$$
    bound suitably framed Seifert surfaces 
    $W^i_{0,1}\subset S^{n+1}\times I$, for details see \S 3 in [Ha,St].
    $\lambda_0$ and $\lambda_1$ are represented by the following 
    framed links in $S^{p_0}$ and $S^{p_1}$ respectively
\begin{align*}
           M_1\sqcup M_2\sqcup\cdots\sqcup M_r
            & =  (S^{p_0}\pitchfork V_1)\sqcup (S^{p_0}\pitchfork V_{1,2})
               \cdots\sqcup (S^{p_0}\pitchfork V_{1,r}), \\ 
           M'_1\sqcup M'_2\sqcup\cdots\sqcup M'_r
            & =  (S^{p_1}\pitchfork V_0)\sqcup (S^{p_1}\pitchfork V_{0,2})
               \cdots\sqcup (S^{p_1}\pitchfork V_{0,r}); 
\end{align*}
    see [Ha,St] again.

       From $V_0\pitchfork V_1$ we get a framed bordism 
    $W_{0,1}\subset S^{n+1}\times I$ between $E^{n+1-p_0}[M_1]$ and 
    $E^{n+1-p_1}[M'_1]$. To see this note that
    $M_1$ lies in a ball $D^{p_0}\subset S^{p_0}$, and we can isotope
    this ball to the standard embedding in $S^{n+1}$ and homotope its
    framing to the standard one, this implies the boundary part
    $M_1$ of $W_{0,1}$ represents $E^{n+1-p_0}[M_1]$ up to sign. 
    For the other boundary
    part it is completely similar. Clearly, the symmetry relation of 
    Kervaire [Ke] is desuspended. 

       We can obviously arrange $W_{0,1}$ so that it is disjoint from 
    $V_{0,i}\times I$ and $V_{1,i}\times I$ for $2\le i\le r$.
    We embed $W^{j_1}_{0,1}$ in the natural way into $S^{n+1}\times I_1$
    and so $W^{j_1}_{0,1}\times I$ into $S^{n+1}\times I\times I_1$.
    Consider now the framed intersection
        $$ Q'_1=W_{0,1}\times I_1\pitchfork W^{j_1}_{0,1}\times I
                \subset S^{n+1}\times I\times I_1. $$
    According to the construction we see directly 
    $Q'_1\subset S^{n+1}\times I\times (0,1)$. It holds in addition
\begin{eqnarray*}
       \partial Q'_1=Z'_1\sqcup Z'_2
         & = & (M_1\times\{0\}\times I_1\pitchfork W^{j_1}_{0,1}\times\{0\})
                \sqcup  \\  \nonumber
         &   & (M'_1\times\{1\}\times I_1\pitchfork           
                W^{j_1}_{0,1}\times\{1\}).  \nonumber
\end{eqnarray*}
    Just as in the $\tau$-construction we isotope $Q'_1$, using the last
    normal vector field $v$ in the framing of $W_{0,1}$, to a framed
    submanifold $Q_1\subset S^{n+1}\times I\times\{0\}$ which lies in
    a small tubular neighbourhood of $W_{0,1}$, $\partial Q_1=Z_1\sqcup Z_2$.
    Up to homotopy of the framing we can assume that $v$ restricts to the
    last vector fields in the framings of $M_1$ and $M'_1$ (considered 
    as submanifolds of $S^{p_0}$ and $S^{p_1}$ respectively). This   
    means $Z_1$ lies in $S^{p_0}$ and is just $\tau(M_1,M_{j_1})$
    up to involution of the framing. In fact the intersection
          $$M_1\times\{0\}\times I_1\pitchfork W^{j_1}_{0,1}\times\{0\}$$
    is just the transversal intersection 
         $$ M_1\times\{0\}\times I_1\pitchfork 
            (S^{p_0}\times\{0\}\times I_1\pitchfork 
            W^{j_1}_{0,1}\times\{0\})   $$
    considered in $S^{p_0}\times\{0\}\times I_1$, in particular
    $S^{p_0}\times\{0\}\times I_1\pitchfork W^{j_1}_{0,1}\times\{0\}$
    is a framed zerobordism of $M_{j_1}$ under the assumption that
    the sublinks $L_0$ and $L_1$ of $L$ are zero $h$-cobordant. 
    The same is true for $Z_2$ (the fact at the beginning of this
    section is used here). So,
    considered in $S^{n+1}$ the framed submanifolds $Z_1$, $Z_2$ 
    represent
       $$\pm E^{n+1-p_0}\tau[M_1,M_{j_1}], \ \ \ \ 
         \pm E^{n+1-p_1}\tau'[M'_1,M'_{j_1}]  $$
    respectively and $Q_1$ is the desired framed bordism. The case 
    $\gamma=[\iota_{j_1},\iota_1]$ follows. Note that the symmetry
    relation of Haefliger and Steer is desuspended one time.

       Let now  $\gamma=[\iota_{j_t},\gamma']$ be as at the beginning
    of this section. Assume inductively that the assertion for 
    $\gamma'$ is true and the corresponding framed bordism $W_{\gamma'}$
    lies in a small tubular neighbourhood of $W_{0,1}$ and considered 
    in $S^{p_0}$ or $S^{p_1}$ its two boundary parts represent 
    $\mu^{\tau}_{\gamma'}(\lambda_0)$ and 
    $\mu^{\tau}_{\gamma'}(\lambda_1)$ respectively. 
    In particularly, this means that
    $W_{\gamma'}$ is disjoint from $V_{0,j_t}\times I$ and 
    $V_{1,j_t}\times I$. The assertion for $\gamma$ follows, if we replace
    $j_1$ and $W_{0,1}$ in the argument above by $j_t$ and $W_{\gamma'}$
    respectively.
    
        II. \ Let $D^{p_0+1}, D^{p_1+1}$ be two disjoint balls in 
    $S^{n+1}$ which are disjoint from all components 
    of the link $L$. Define
\begin{align*}
          L'_0 & = \partial D^{p_0+1}\sqcup S^{p_1}\sqcup S^{p_2}\sqcup
                     \cdots \sqcup S^{p_r}  \\   
          L'_1 & = S^{p_0}\sqcup\partial D^{p_1+1}\sqcup S^{p_2}\sqcup   
                     \cdots \sqcup S^{p_r}  \\
      L'_{0,1} & = \partial D^{p_0+1}\sqcup\partial D^{p_1+1}\sqcup 
                     S^{p_2}\sqcup\cdots \sqcup S^{p_r}.  
\end{align*}
    If the condition in part I is not satisfied, namely if the sublinks
    $L_0$ and $L_1$ of $L$ are not zero $h$-cobordant, then consider the 
    connected sum
            $$ L'=L-L'_0-L'_1+L'_{0,1}.  $$
    $L'$ satisfies clearly the just mentioned condition. For example, by    
    forgetting the $p_0$-dimensional sphere in $L'$ we obtain
\begin{align*}
           & (S^{p_1}\sqcup S^{p_2}\sqcup\cdots \sqcup S^{p_r})-
             (S^{p_1}\sqcup S^{p_2}\sqcup\cdots \sqcup S^{p_r})-  \\
           & (\partial D^{p_1+1}\sqcup S^{p_2}\sqcup\cdots \sqcup S^{p_r})+ 
             (\partial D^{p_1+1}\sqcup S^{p_2}\sqcup\cdots \sqcup S^{p_r})
\end{align*}
    which is evidently zero $h$-cobordant. Denote by 
        $$ \lambda_0, \lambda_0^{(1)}, \lambda_0^{(2)}, \lambda_0^{(3)},
           \lambda_0^{(4)}\in\pi_{p_0}(S^{k_1}\vee\cdots\vee S^{k_r}) $$
    the elements given by $ L, L', L'_0, L'_1$ and $L'_{0,1}$ respectively.
    It holds $\lambda_0^{(2)}=\lambda_0^{(4)}=0$ and therefore 
    $\lambda_0^{(1)}=\lambda_0-\lambda_0^{(3)}$. Because the first 
    component of the link representing $\lambda_0^{(3)}$ is the empty,
    we see $\mu_{\gamma}^{\tau}(\lambda_0^{(1)})=
    \mu_{\gamma}^{\tau}(\lambda_0)$. $\mu_{\gamma}^{\tau}(\lambda_1^{(1)})=
    \mu_{\gamma}^{\tau}(\lambda_1)$ follows by the same token. We finish
    the proof by using part I.
\end{proof}

      We will obtain more symmetry relations if we replace the pair
    $(0,1)$ by any $(i,j)$ with $0\le i\not=j\le r$. We do not know
    the exact relationship between our symmetry relations and those 
    of Turaev in [Tu] and those of Koschorke in [Ko 3]. 
    We presume that our relations
    can in general not be desuspended, because the framed manifolds
    like $W_{0,1}$ and $W^i_{0,1}$, which we have used, take their 
    place very naturally in the sphere $S^{n+1}$.

\vskip 1cm
\section{the $\tau$-reduction }
\vskip 3pt

        In this section we define first the $\tau$-reductions by using 
    $\tau$-constructions and then construct all the Hilton homomorphisms
    geometrically by means of $\tau$-reductions.

      Let $\gamma$ be a Whitehead product in $\iota_1,\cdots, \iota_r$
    and $1\le i<j\le r$. If we replace all $[\iota_i,\iota_j]$ and
    $[\iota_j,\iota_i]$ in $\gamma$ by $\iota_{r+1}$, then we get a 
    new Whitehead product $\tau^S_{ji}(\gamma)$ in 
    $\iota_1,\cdots, \iota_r, \iota_{r+1}$. We call $\tau^S_{ji}$ a
    symbolic reduction. Note that $\tau^S_{ji}(\gamma)$ is generally
    not a basic Whitehead product even if $\gamma$ is.

       We will construct by geometrical means a homomorphism $\tau^R_{ji}$
    such that the following diagram commutes for some Whitehead
    products $\gamma$. We call $\tau^R_{ji}$ a $\tau$-reduction.
    In the diagram $k_{r+1}=k_i+k_j-1$.
\[
\begin{diagram}
        \node{FL_X^{k_1,\cdots,k_r}}
        \arrow{e,t}{\tau_{ji}^R}
        \node{FL_X^{k_1,\cdots,k_r,k_{r+1}}} 
\\
        \node{FL_X^{q(\gamma)+1}}
        \arrow{n,l}{\gamma_*}
        \arrow{ne,r}{(\tau_{ji}^S(\gamma))_*}
\end{diagram}
\]

       Fix $m$, the dimension of $X$, and the codimensions 
    $ k_1, \cdots, k_r$. All these numbers are supposed to be $\ge 2$. 
    Then there 
    is a $w_0$ such that $FL_X^{q(\gamma)+1}=0$ holds for all Whitehead
    products $\gamma$ of weight greater than $w_0$. We define for
    $2\le w\le w_0$ and $1\le i<j\le r$ a homomorphism $\tau_{ji}^w$ as
    follows. Consider a framed $(k_1, \cdots, k_r, k_{r+1})$-link
    $M_1\sqcup M_2\sqcup\cdots\sqcup M_r\sqcup M_{r+1}$ and define
\begin{align*}
          Z'_w &= \tau(M_i,\cdots\tau(M_i,\tau(M_i,M_i^{sh}))\cdots), \\
          Z_w  &= \tau(M_j,Z'_w),  
\end{align*}
    where $M_i^{sh}$ is a small shift of $M_i$ along the framing and
    $(w-2)$ $\tau$-constructions are used to get $Z'_w$. If $w=2$ let
    $Z'_w=M_i$. Let 
    $\gamma_w=[\iota_i,[\iota_i,\cdots[\iota_i,\iota_j]\cdots]]$ 
    be of weight $w$ and let $L_w, L_w^S$ be the framed links in 
    $\Bbb R^{q(\gamma_w)+1}$ representing $\gamma_w$ and 
    $\tau_{ji}^S(\gamma_w)$ respectively. The framed submanifold 
    $Z_w$ is of codimension $q(\gamma_w)+1$, so we can embed the sum
    $-(L_w\sqcup\phi)+L_w^S$ fibre-wise in a small tubular neighbourhood
    of $Z_w$ to get a new framed link
\begin{eqnarray}
        M'_1\sqcup M'_2\sqcup\cdots\sqcup M'_r\sqcup M'_{r+1}
           =Z_w\times (-(L_w\sqcup\phi)+L_w^S).
\end{eqnarray}
    Note, by some suitable conventions of the framings involved we can 
    guarantee that for the link representing $\gamma_w$ it holds 
    strictly $Z_w=+1$ (that $Z_w$ is a single point is proved in Lemma 4.3), 
    from now on we assume this has been done. We define now
\begin{eqnarray*}
       \tau_{ji}^w(M_1\sqcup\cdots\sqcup M_r\sqcup M_{r+1})
          &=& \bar M_1\sqcup\cdots\sqcup \bar M_r\sqcup\bar M_{r+1} \\
          &=& (M_1\sqcup M'_1)\sqcup\cdots
              \sqcup (M_{r+1}\sqcup M'_{r+1}). 
\end{eqnarray*}
\begin{Lemma} \  The assignment $\tau_{ji}^w$ above gives a well
     defined homomorphism
        $$\tau_{ji}^w:FL_X^{k_1,\cdots,k_r,k_{r+1}}\longrightarrow
                      FL_X^{k_1,\cdots,k_r,k_{r+1}}.    $$
\end{Lemma}
\begin{proof} \ The $\tau$-constructions, fibre-wise embeddings and fusion
     of components are or induce well defined homomorphisms.
\end{proof}
\begin{Definition} \  Let $incl_*:FL_X^{k_1,\cdots,k_r}
     \longrightarrow FL_X^{k_1,\cdots,k_r,k_{r+1}}$ be the inclusion,
     where $k_{r+1}=k_i+k_j-1$. We define 
       $$ \tau_{ji}^R=\tau_{ji}^2\circ\tau_{ji}^3\circ\cdots
          \circ\tau_{ji}^{w_0}\circ incl_*$$  
     and call it a $\tau$-reduction.                              
\end{Definition}

        We observe the following: if $M_i$ is framed zerobordant, then
     the link $M_i\sqcup M_i^{sh}$ is framed zerobordant and therefore
     every $Z'_w$ defined as above is framed zerobordant.
\begin{Lemma} \ Let $\gamma_w$, $L_w$ and $L_w^S$ be as above. It holds
     $\tau_{ji}^R(L_w)=L_w^S$.
\end{Lemma}
\begin{proof} \ The link $L_w$ is given by (3) and (4) 
     in the proof of Theorem 2.3. Part I of this proof is heavily 
     based on the following observations:
\begin{quote}  
          (i) \ Every $N_{i,\lambda}\subset M_i$ lies in a 
          $q(\gamma_w)$-dimensional subspace of $\Bbb R^{q(\gamma_w)+1}$
          and has therefore a (constant) vector in $\Bbb R^{q(\gamma_w)+1}$
          as the last normal vector field in its framing (at least up to
          homotopy of the framing), $1\le \lambda\le w-1$.
\newline
          (ii) \ Every sub-product $P_{i,\lambda}^{sub}\subset N_{i,\lambda}$
          containing the factor $S^{(\lambda-1)(k_i-1)+k_j-1}$ bounds
          a suitably framed Seifert surface $F_{i,\lambda}^{sub}$.
          If $\lambda_1<\lambda_2$ then we have 
          $F_{i,\lambda_1}^{sub}\pitchfork P_{i,\lambda_2}^{sub}=\phi$.
          If $\lambda_1=\lambda_2$ we may shift one of the two sub-products 
          slightly along the framing and see that the same is true 
          according to (i).
\end{quote}
       
        I. \ We prove first the following assertion by induction: evaluated
     on $L_w$ it holds $Z'_{\bar w}=\phi$ for $\bar w>w$ and ($\{pt\}$
     is a set consisting of a single point)
            $$ Z'_w=\{pt\}\times\cdots\times\{pt\}\times S^{k_j-1},$$
     which is a small shift of the obvious submanifold in $N_{i,1}$.
     For $w=2$ it is trivial, so assume $w\ge 3$. Let
     $\hat M_k$ be the components of the link representing $\gamma_{w-1}$.
     According to (3) and (4) in the proof of Theorem 2.3,  
     $\hat M_i=\hat N_{i,1}\sqcup\cdots\sqcup\hat N_{i,w-2}$ and
\begin{eqnarray*}
           M_i & = & S^{k_i-1}_{w-1}\times 
                     \hat M_i\sqcup S^{(w-2)(k_i-1)+k_j-1} \\
               & = & S^{k_i-1}_{w-1}\times \hat N_{i,1}
                    \sqcup\cdots\sqcup
                    S^{k_i-1}_{w-1}\times \hat N_{i,w-2}
                    \sqcup S^{(w-2)(k_i-1)+k_j-1}  \\
               & = &  :N_{i,1}\sqcup\cdots\sqcup N_{i,w-1}, \\
          M_j  & = & S^{k_i-1}_{w-1}\times \hat M_j. 
\end{eqnarray*} 
     All other components are the empty. Using the observations we get
\begin{align*}
           Z'_3 &= \tau(M_i,M_i^{sh})     \\
                &= S^{k_i-1}\times \hat Z'_3\sqcup
                    (\sqcup_{\lambda_1=1}^{w-2}
                    \tau(N_{i,\lambda_1},N_{i,w-1}^{sh})),
\end{align*}
     where $\hat Z'_3=\tau(\hat M_i,\hat M_i^{sh})$ and $N_{i,w-1}^{sh}$
     is a small shift of $N_{i,w-1}$. We use the observation again
     and get
          $$ Z'_4=S^{k_i-1}\times \hat Z'_4\sqcup
                    (\sqcup_{\lambda_2=1}^{w-3}
                    \sqcup_{\lambda_1>\lambda_2}^{w-2}
                    \tau(N_{i,\lambda_2},
                    \tau(N_{i,\lambda_1},N_{i,w-1}^{sh}))).  $$
     Just repeat this until we get $Z'_w$. 
     $S_{w-1}^{k_i-1}\times \hat Z'_w=\phi$ is obvious (using the induction
     assumption). Denote by $\Lambda$ the condition
        $$ 1\le \lambda_{w-2}<\lambda_{w-3}<\cdots<
           \lambda_2<\lambda_1\le w-2.$$
     The other part of $Z'_w$ is given by
\begin{eqnarray*}         
          &   &  \sqcup_{\Lambda}\tau(N_{i,\lambda_{w-2}},\tau(\cdots,
                \tau(N_{i,\lambda_2},
                \tau(N_{i,\lambda_1},N_{i,w-1}^{sh}))\cdots))  \\       
           & =&  \tau(N_{i,1},\tau(\cdots,
                \tau(N_{i,w-3},\tau(N_{i,w-2},N_{i,w-1}^{sh}))\cdots)) \\
           & =& \{pt\}\times\cdots\times\{pt\}\times S^{k_j-1}.
\end{eqnarray*}
     For $\bar w>w$ we see $Z'_{\bar w}=\phi$ immediately.

         II. \ From part I we obtain $Z_w=\{pt\}$ with positive sign, 
     and $Z_{\bar w}=\phi$ if $\bar w>w$. This means 
         $$ \tau_{ji}^{w+1}\circ\cdots\circ\tau_{ji}^{w_0}\circ incl_*
            (L_w)=L_w\sqcup\phi. $$
     From the definition of $\tau_{ji}^w$ we also have
         $$ \tau_{ji}^w(L_w\sqcup\phi)=(L_w\sqcup\phi)+
            \{pt\}\times (-(L_w\sqcup\phi)+L_w^S)=L_w^S.  $$
     Because the $j$-th component of $L_w^S$ is empty it follows
     $\tau_{ji}^{2}\circ\cdots\circ\tau_{ji}^{w-1}(L_w^S)=L_w^S$.
     This shows $\tau_{ji}^R(L_w)=L_w^S$.
\end{proof}

        We hope that the background of the definition of $\tau_{ji}^R$
     is more or less presented in the proof of this lemma. Recall 
     formula (5). We use the nagative part $-(L_w\sqcup\phi)$ to eliminate
     what troubles us and use the part $L_w^S$ to get what we desire. We show
     next that the $\tau$-reductions fit in the commutative diagram
     at the beginning of this section for some Whitehead products.
\begin{Lemma} \ Let $\Gamma$ be a system of basic Whitehead products
     in $\iota_1<\iota_2<\cdots<\iota_r$. For all basic Whitehead products
     $\gamma\in\Gamma$ it holds   
     $\tau_{j,1}^R\circ\gamma_*=(\tau_{j,1}^S(\gamma))_*$.
\end{Lemma}
\begin{proof} \ If $\gamma$ is of weight $1$ then the statement is trivial.
     Let the weight of $\gamma$ be greater than $1$ and let
     $L_{\gamma}=M_1\sqcup M_2\sqcup\cdots\sqcup M_r$ be the framed link
     representing $\gamma$. Because every component $M_i$ is framed 
     zerobordant we need only to show     
     $\tau_{j,1}^R(L_{\gamma})=L_{\gamma}^S$ according to Lemma 2.2, where
     $L_{\gamma}^S$ represents $\tau_{j,1}^S(\gamma)$. If $\gamma$ is of
     weight $2$ the statement follows easily. Let $\gamma=[\alpha,\beta]$
     be of weight $\ge 3$. Then formula (1) in the proof of Theorem 2.3
     holds. According to the definition $\beta$ has weight at least $2$,
     therefore all components $M_k(\beta)$ of the framed link $L_{\beta}$
     representing $\beta$ are framed zerobordant. If $\alpha\not=\iota_1$
     then the first component $M_1(\alpha)$ of the link $L_{\alpha}$
     representing $\alpha$ is also framed zerobordant. Using Lemma 2.2
     we get the following formula similar to (2) in the proof of Theorem 2.3
       $$ \tau_{j,1}^R(L_{\gamma})=S^{l_1-1}\times\tau_{j,1}^R(L_{\beta})
             \sqcup S^{l_2-1}\times\tau_{j,1}^R(L_{\alpha}).  $$
     Under the inductive assumption for $\alpha$ and $\beta$ the assertion
     follows from this formula.

        If $\alpha=\iota_1$ then the only possiblity is
     $\gamma=[\iota_1,[\iota_1,\cdots[\iota_1,\iota_{j'}]\cdots]]$
     according to the definition of basic Whitehead products.
     If $\iota_{j'}\not=\iota_j$ then $\iota_j$ does not appear in
     $\gamma$, for $\iota_1<\iota_j$. This means the $j$-th component
     $M_j$ is the empty and the statement follows. If  
     $\iota_{j'}=\iota_j$ then Lemma 4.3 applies.
\end{proof}

        We define now an ordered sequence $T^S$ of symbolic reductions
    as follows. Let
       $$\iota_1<\cdots <\iota_r<\gamma_1<\cdots
         <\gamma_n<\gamma_{n+1}<\cdots $$
    be the sequence of basic Whitehead products in $\Gamma$. 
    If $\gamma_1=[\iota_{i_1},\iota_{j_1}]$ then the first reduction in
    $T^S$ is $\tau_{j_1,i_1}^S$ determined by $\gamma_1$. After this 
    reduction we get
        $$\iota_1<\cdots <\iota_r<\iota_{r+1}<\gamma_2^1<\cdots<
          \gamma_n^1<\gamma_{n+1}^1<\cdots $$
    Now $\gamma_2^1=[\iota_{i_2},\iota_{j_2}]$ is a Whitehead product
    in $\iota_1, \cdots, \iota_r, \iota_{r+1}$ of weight $2$ from which
    we obtain the second reduction $\tau_{j_2,i_2}^S$ in $T^S$. After
    doing this $n$-times one gets
        $$\iota_1<\cdots <\iota_r<\cdots<\iota_{r+n}<\gamma_{n+1}^n
          <\gamma_{n+2}^n<\cdots <\gamma_{m}^n<\gamma_{m+1}^n\cdots $$
    It is not difficult to show that for all $k\ge 1$ the Whitehead
    products $\gamma_{n+k}^n$ are different and are of weight $\ge 2$,
    and $\gamma_{n+1}^n=[\iota_{i_{n+1}},\iota_{j_{n+1}}]$
    is of weight $2$. So we define the $(n+1)$-th reduction to be
    $\tau_{j_{n+1},i_{n+1}}^S$. Defined in this way $T^S$ reduces
    the original sequence to the following 
         $$\iota_1<\cdots <\iota_r<\cdots<\iota_{r+n}<\iota_{r+n+1}\cdots $$
    Note that if the numbers $m, k_1, \cdots, k_r$ are fixed, then
    the sequence of basic Whitehead products $\gamma_i$ 
    and the sequence $T^S$ of reductions are both finite.
\begin{Definition} \ Define $T^R$ to be the sequence of 
    $\tau$-reductions determined by $T^S$. We denote by $T^S_n$ 
    and $T^R_n$ the first $n$ reductions in $T^S$ and $T^R$ respectively.
\end{Definition}
\begin{Theorem} \ Let $\Gamma$ be a system of basic Whitehead
    products in $\iota_1<\cdots <\iota_r$, such that the conditions
    $\alpha_1<\alpha_2$ and $[\alpha_1,\beta], [\alpha_2,\beta]\in\Gamma$
    together imply $[\alpha_1,\beta]<[\alpha_2,\beta]$. Then
    the diagram
\[
\begin{diagram}
        \node{FL_X^{k_1,\cdots,k_r}}
        \arrow{e,t}{T^R}
        \node{FL_X^{k_1,\cdots,k_r,\cdots,k_{r+n(\gamma)},\cdots}} 
        \arrow{s,r}{p_{\gamma}}
\\
        \node{FL_X^{q(\gamma)+1}}
        \arrow{e,b}{id}
        \arrow{n,l}{\gamma_*}
        \arrow{ne,r}{(\iota_{r+n(\gamma)})_*}
        \node{FL_X^{q(\gamma)+1}}
\end{diagram}
\]
    commutes for all $\gamma\in\Gamma$, where $k_{r+n(\gamma)}=q(\gamma)+1$,
    $p_{\gamma}$ is the obvious projection, $\iota_{r+n(\gamma)}=T^S(\gamma)$
    and we have assumed $\gamma$ is the $(r+n(\gamma))$-th basic 
    Whitehead product in $\Gamma$. In particular, 
    $\Delta_{\gamma}=p_{\gamma}\circ T^R$ is exactly the Hilton
    homomorphism $H_{\gamma}$.
\end{Theorem}
\begin{proof} \ If the first statement is true then one can easily check
    that $\Delta_{\gamma}$ satisfies the properties (a) and (b) in
    \S 1. For example
          $$  \Delta_{\gamma}\circ\gamma_*
             =p_{\gamma}\circ T^R\circ\gamma_*
             =p_{\gamma}\circ (\iota_{r+n(\gamma)})_*=id.  $$
    (b) follows easily. This shows $H_{\gamma}=\Delta_{\gamma}$.

        For the first statement we need to show 
    $T^R_k\circ\gamma_*=(T^S_k(\gamma))_*$ for any $k$. 
    Because $\iota_1$ is the
    first basic Whitehead product in $\Gamma$ the first reduction is
    $\tau_{j_1,1}^R$. So the case $k=1$ follows from Lemma 4.4. Assume
    inductively that the statement is true for all natural numbers
    $\le k$. We will prove the case $k+1$ by induction on the weight
    $w$ of $T^S_k(\gamma)$.

        Let $\tau_{ji}^R$ be the $(k+1)$-th reduction. If $w=1$ the
    assertion is trivial. Let $w\ge 2$ and denote by $L_{\gamma}^k$ the 
    framed link representing $T^S_k(\gamma)$. The components of this 
    link are framed zerobordant. By Lemma 2.2 we just need to show 
    $\tau_{ji}^R(L_{\gamma}^k)=L_{\gamma}^{k+1}$. If $w=2$ this is not
    difficult to see. Let then $T_k^S(\gamma)=[T_k^S(\alpha),T_k^S(\beta)]$
    to be of weight $\ge 3$, and $\gamma=[\alpha,\beta]\in\Gamma$. Then
    we have
       $$ M_l(\gamma)=S^{l_1-1}\times M_l(\beta)
                      \sqcup S^{l_2-1}\times M_l(\alpha), $$
    where $ M_l(\alpha), M_l(\beta), M_l(\gamma)$ are components of
    the links $L_{\alpha}^k, L_{\beta}^k, L_{\gamma}^k$ representing
    $T^S_k(\alpha)$, $T^S_k(\beta)$ and $T^S_k(\gamma)$ 
    respectively. According to the
    definition of $T^S$ and the basic Whitehead products we know 
    $T^S_k(\beta)$ is at leat of weight $2$ and therefore the components
    of the corresponding link are framed zerobordant. If      
    $T^S_k(\alpha)\not=\iota_i$ then $M_i(\alpha)$ is also framed 
    zerobordant. Using Lemma 2.2 again we obtain 
         $$ \tau_{ji}^R(L^k_{\gamma})=S^{l_1-1}\times\tau_{ji}^R(L^k_{\beta})
             \sqcup S^{l_2-1}\times\tau_{ji}^R(L^k_{\alpha}).  $$
    The statement for $\gamma$ now follows from this formula under the 
    inductive assumption for $\alpha$ and $\beta$. If
    $T^S_k(\alpha)=\iota_i$ then
         $$  T^S_k(\gamma)=[\iota_{i},[\iota_{i_1},\cdots
             [\iota_{i_t},\iota_{j'}]\cdots ]],  $$
    with $\iota_i\ge \iota_{i_1}\ge\cdots\ge\iota_{i_t}<\iota_{j'}$
    and $\iota_i<\iota_j$, according to the construction of $T^S$.
    So if $\iota_{j'}\not=\iota_j$ then $\iota_j$ does not appear in
    $T_k^S(\gamma)$ and the assertion follows trivially. Let
    $\iota_{j'}=\iota_j$ and assume $\iota_i>\iota_{i_t}$. Denote
    the original basic Whitehead products of $[\iota_{i_t},\iota_j]$
    and $[\iota_{i},\iota_j]$ by $\gamma_1, \gamma_2$ respectively.
    Then $\gamma_1<\gamma_2$, according to the condition on
    $\Gamma$ in the theorem. This means
    $\tau^R_{j,i_t}$ is the $k'$-th reduction with $k'<k+1$. But
    after this $k'$-th reduction there is no appearance of
    $[\iota_{i_t},\iota_j]$. So $\iota_i=\iota_{i_1}=\cdots =\iota_{i_t}$
    is the only possibility. The statement follows now from Lemma 4.3.
\end{proof}
 
       The restriction in Theorem 4.6 on the system of basic Whitehead  
    products is not necessary, but the sequences $T^S$ and $T^R$ should
    be adjusted as follows. Let $\Gamma$ be any system of basic Whitehead
    products in $\iota_1<\iota_2\cdots <\iota_r$. We give the elements 
    of $\Gamma$ a new order $\prec$ which respects the weights and has
    the following
    property: if $\alpha_1, \alpha_2, [\alpha_1,\beta], [\alpha_2,\beta]$
    are basic Whitehead products in $\Gamma$ and $\alpha_1<\alpha_2$, 
    then $[\alpha_1,\beta]\prec[\alpha_2,\beta]$. 
    It is easily seen that such an order does exist.
    Using this new order of $\Gamma$ we get a sequence $T^S(\Gamma)$
    of symbolic reductions and the corresponding sequence $T^R(\Gamma)$
    of $\tau$-reductions. If $\gamma$ is the $(r+n(\gamma))$-th
    and $(r+n'(\gamma))$-th basic Whitehead product in $\Gamma$ with
    respect to the old order $<$ and the new order 
    $\prec$ respectively, then the $(r+n'(\gamma))$-th reduction 
    in $T^S(\Gamma)$ is the one after which 
    $\gamma$ is just reduced to $\iota_{r+n(\gamma)}$ 
    (not $\iota_{r+n'(\gamma)}$). Note that if $\gamma$ and $\gamma'$ 
    are basic Whitehead products of the same weight and if
    the first $k$ reductions $T^S_k(\Gamma)$ already reduce $\gamma$ to
    weight $1$, then the weight of $T^S_k(\Gamma)(\gamma')$ must $\le 2$.
\begin{Proposition} \ The result in Theorem 4.6 still holds
    for any system $\Gamma$ of basic Whitehead products in
    $\iota_1<\iota_2\cdots <\iota_r$, if we replace
    the sequences $T^S$, $T^R$ there by the sequences $T^S(\Gamma)$   
    and $T^R(\Gamma)$.
\end{Proposition}
\begin{proof} \ Except following changes the proof remains the same.

        Let $\tau_{ji}^R$ be the $(k+1)$-th reduction and  
    $T_k^S(\Gamma)(\gamma)=[T_k^S(\Gamma)(\alpha),T_k^S(\Gamma)(\beta)]$
    be of weight $\ge 3$, $\gamma=[\alpha,\beta]\in\Gamma$.
    We see $T_k^S(\Gamma)(\beta)$ may be of weight $1$, in this case
    $T_k^S(\Gamma)(\alpha)$ must be of weight $2$ according to the
    definitions of basic Whtehead products and the sequence of 
    symbolic reductions. So  
    $T_k^S(\Gamma)(\gamma)=[[\iota_{i_1},\iota_{i_2}],\iota_{i_3}]$.
    Let $\iota_{i_3}=\iota_i$, then no one of $\iota_{i_1},\iota_{i_2}$
    can be $\iota_j$, because, if $\alpha_1,\alpha_2,\gamma'$ are the
    original basic Whitehead products of $\iota_{i_1},\iota_{i_2}$
    and $\iota_j$, then we have $w(\alpha_1)\le w(\alpha_2)<w(\alpha)\le   
    w(\beta)\le w(\gamma')$ according to the definition of the 
    basic Whitehead products. So the assertion follows trivially.
    Other cases (where $\iota_{i_3}\not=\iota_i$) can be easily checked.
     
       If the weights of both  $T_k^S(\Gamma)(\alpha)$
    and $T_k^S(\Gamma)(\beta)$ are at least $2$, or if 
    $T_k^S(\Gamma)(\alpha)$ is of weight $1$ but is different
    from $\iota_i$ and $T_k^S(\Gamma)(\beta)$ is of weight at least $2$,
    then the statement follows by induction as in the proof of Theorem 4.6.

       So assume the weight of $T_k^S(\Gamma)(\beta)$ is at least $2$
    and $T_k^S(\Gamma)(\alpha)=\iota_i$, then
       $$ T_k^S(\Gamma)(\gamma)=[\iota_i,[T_k^S(\Gamma)(\beta_1),
          T_k^S(\Gamma)(\beta_2)]].  $$
    Because $\alpha\ge \beta_1$, the weight of the original basic 
    Whitehead product of $[\iota_i,\iota_j]$ is clearly greater than 
    the weight of $\beta_1$; according to the definitions of $\prec$
    and $T^S(\Gamma)$, a basic Whitehead product can only be symbolically 
    reduced to weight $1$ when all the basic Whitehead products of smaller 
    weight are already reduced to weight $1$. This means
    $T_k^S(\Gamma)(\beta_1)=\iota_{i_1}$ must be of weight $1$. Therefore 
       $$ T_k^S(\Gamma)(\gamma)=[\iota_i,[\iota_{i_1},\cdots,
         [\iota_{i_t},\iota_{j'}]\cdots]]. $$ 
    
       Now let $\alpha,\beta'',\alpha_1,\cdots,\alpha_t,\beta'$ be the
    original basic Whitehead products corresponding to
    $\iota_i,\iota_j,\iota_{i_1},\cdots,\iota_{i_t},\iota_{j'}$. Because
       $$ \beta''>\alpha\ge \alpha_1\ge\cdots\ge\alpha_t<\beta', $$
    if $\iota_{j'}\not=\iota_j$ then no one of 
    $\iota_i,\iota_{i_1},\cdots,\iota_{i_t}$ can be $\iota_j$,
    the statement follows trivially. 
  
       Let now $\iota_{j'}=\iota_j$. If $\iota_{i_t}\not=\iota_i$,
    then $[\alpha_t,\beta]\prec [\alpha,\beta]$ with respect to the
    new order, because $\alpha_t<\alpha$ with respect to the original
    order. By the definition of $T^S(\Gamma)$ this means
    $\tau^S_{j, i_t}$ is the $k'$-th reduction in $T^S(\Gamma)$
    with $k'<k+1$. But after this $k'$-th reduction 
    $[\iota_{i_t},\iota_j]$ can not appear in $T_k^S(\Gamma)(\gamma)$ for
    any $\gamma$. So the only possibility is 
        $$ T_k^S(\Gamma)(\gamma)=[\iota_i,[\iota_{i},\cdots,
         [\iota_{i},\iota_{j}]\cdots]]. $$
    The assertion follows now from Lemma 4.3. 
\end{proof}

       Note that, if $T^S_n(\Gamma)$ already reduces $\gamma$ to
    weight $1$, then $H_{\gamma}=\Delta^n_{\gamma}
    =p_{\gamma}\circ T^R_n(\Gamma)$; and if $\tau^R_{ji}$ is the $(k+1)$-th
    reduction and at least one of $\iota_i$ and $\iota_j$ does not
    appear in $T^S_k(\Gamma)(\gamma)$, we can eliminate $\tau^R_{ji}$
    from $T^R_n(\Gamma)$. Therefore if the weight of $\gamma$ is $w$, 
    we need exactly $w-1$ $\tau$-reductions to get $H_{\gamma}$.
  
\begin{Example} \ Let $r=2$ and consider 
        $$ \iota_1<\iota_2<\gamma_1<\gamma_2<\gamma_3<\gamma_4
           <\gamma_5<\gamma_6<\cdots $$ 
    where $\gamma_1=[\iota_1,\iota_2]$, $\gamma_2=[\iota_1,[\iota_1,\iota_2]]$,
    $\gamma_3=[\iota_2,[\iota_1,\iota_2]]$,  
    $\gamma_4=[\iota_1,[\iota_1,[\iota_1,\iota_2]]]$. Then the 
    reductions $T^S_4=(\tau^S_{2,1}, \tau^S_{3,1},\tau^S_{3,2},\tau^S_{4,1})$
    reduce the sequence above to
        $$ \iota_1<\iota_2<\iota_3<\iota_4
           <\iota_5<\iota_6<\gamma^4_5<\gamma^4_6\cdots $$ 
    So, the Hilton homomorphisms corresponding to $\gamma_1$, $\gamma_2$,
    $\gamma_3$, $\gamma_4$ are given by
        $$ p_3\circ\tau^R_{2,1}, \hskip 1cm
           p_4\circ\tau^R_{3,1}\circ\tau^R_{2,1},  $$
        $$ p_5\circ\tau^R_{3,2}\circ\tau^R_{2,1},\hskip  1cm 
           p_6\circ\tau^R_{4,1}
                \circ\tau^R_{3,1}\circ\tau^R_{2,1},          $$
    where $p_3, p_4, p_5$ and $p_6$ are the obvious projections.
\end{Example}

       Let $B_i$ be a closed connected smooth manifold and $\xi_i$ be a
    differential vector bundle over $B_i$, $i=1, 2, \cdots, r$.
    We may consider links $M_1\sqcup \cdots \sqcup M_r$ in $X\times \Bbb R$    
    such that the normal bundle of $M_i$ is classified by a bundle
    map into $\xi_i\oplus \varepsilon$, where $\varepsilon$ is the
    trivial line bundle. We call such links 
    $(\xi_1, \cdots, \xi_r)$-links. For the bordism group of 
    such links we have the Hilton-Milnor splitting.
    Note that to perform our 
    $\tau$-construction we need only one normal vector field,
    so $\tau$-construction can easily be generalized to
    $(\xi_1, \cdots, \xi_r)$-links. Some basic things concerning 
    this have been done in the author's dissertation [Wa]. 
    We presume there are no essential difficulties to generalize 
    the discussions in this section to the mentioned case above. 
  
\vskip 1.5cm
{\it
Address

Department of Applied Mathematics

Shanxi University of Finance and Economics

030006 Taiyuan, Shanxi, PR China
\begin{verbatim}  e-mail: dkwjh@public.ty.sx.cn \end{verbatim}

\vskip 0.3cm
Currently:

Department of Mathematics

Bar-Ilan University

52900 Ramat Gan, Israel
\begin{verbatim}  e-mail: wangji@macs.biu.ac.il  \end{verbatim}  }

\vskip 20pt
%\Refs
%\widestnumber\key{Ko,Sa 2}
%\vskip 0.2 true cm

\end{document}